\documentclass[twocolumn]{article}

\usepackage{hyperref}
\usepackage{mathtools}
\usepackage{amsfonts}
\usepackage{booktabs}
\usepackage{subfigure}
\usepackage[ruled,vlined,linesnumbered]{algorithm2e}
\usepackage{amsmath}
\usepackage{graphicx}
\usepackage{accents}
\usepackage{listings}
\usepackage{xcolor}
\usepackage[numbers]{natbib}
\newcommand{\dfun}[0]{\Delta^{\ast}}
\newcommand{\under}[1]{\underaccent{\bar}{#1}}

\newcommand{\feas}[0]{\mathit{feas}}
\SetKwRepeat{Do}{do}{while}

\lstdefinestyle{mystyle}{
    backgroundcolor=\color{lightgray!10},   
    basicstyle=\ttfamily\tiny,
    breaklines=true,
    frame=single
}
\title{Certified Inductive Synthesis for Online Mixed-Integer Optimization}

\author{
    Marco Zamponi$^1$, Emilio Incerto$^1$, Daniele Masti$^2$, Mirco Tribastone$^1$ \thanks{© \{Authors | ACM\}, 2025. This is the author's version of the work. It is posted here for your personal use. Not for redistribution. The definitive version was accepted for publication in the proceedings of the ACM/IEEE 16th International Conference on Cyber-Physical Systems (ICCPS '25) \url{https://doi.org/10.1145/3716550.3722014}}\\
    $^1$IMT School for Advanced Studies, Lucca, Italy \\
    $^2$GSSI Gran Sasso Science Institute, L'Aquila, Italy \\
    \texttt{\{marco.zamponi, emilio.incerto, mirco.tribastone\}@imtlucca.it} \\
    \texttt{daniele.masti@gssi.it}
}

\date{}

\begin{document}

\maketitle

\begin{abstract}
    
In fields such as autonomous and safety-critical systems, online optimization plays a crucial role in control and decision-making processes, often requiring the integration of continuous and discrete variables. These tasks are frequently modeled as mixed-integer programming (MIP) problems, where feedback data are incorporated as parameters. However, solving MIPs within strict time constraints is challenging due to their $\mathcal{NP}$-complete nature. A promising solution to this challenge involves leveraging the largely invariant structure of these problems to perform most computations offline, thus enabling efficient online solving even on platforms with limited hardware capabilities. In this paper we present a novel implementation of this strategy that uses counterexample-guided inductive synthesis to split the MIP solution process into two stages. In the offline phase, we construct a mapping that provides feasible assignments for binary variables based on parameter values within a specified range. In the online phase, we solve the remaining continuous part of the problem by fixing the binary variables to the values predicted by this mapping. Our numerical evaluation demonstrates the efficiency and solution quality of this approach compared to standard mixed-integer solvers, highlighting its potential for real-time applications in resource-constrained environments.

\end{abstract}

\bigskip\hrule
\footnotesize
\noindent Please cite this version of the paper:

\noindent Marco Zamponi, Emilio Incerto, Daniele Masti, and Mirco Tribastone, "Certified Inductive Synthesis for Online Mixed-Integer Optimization," in Proceedings of the ACM/IEEE 16th International Conference on Cyber-Physical Systems (ICCPS), 2025, doi: 10.1145/3716550.3722014.

You can use the following BibTeX entry:
\begin{lstlisting}
@inproceedings{zamponi2025certified,
  title={Certified Inductive Synthesis for Online Mixed-Integer Optimization},
  author={Zamponi, Marco and Incerto, Emilio and Masti, Daniele and Tribastone, Mirco},
  booktitle={2025 ACM/IEEE 16th International Conference on Cyber-Physical Systems (ICCPS)},
  year={2025},
  organization={ACM},
  doi={10.1145/3716550.3722014}
}
\end{lstlisting}
\normalsize
\hrule

\section{Introduction}
\label{sec:introduction}

Mathematical programming is a powerful tool to represent decision-making problems. Mixed-integer programming (MIP), in particular,
has several applications in a wide spectrum of domains including job scheduling~\cite{ku2016mixed, pop2016design}, cyber-physical systems~\cite{bemporad1999control}, and neural network verification~\cite{tjeng2017evaluating}. 
Unfortunately, finding an optimal solution for MIP problems is known to be $\mathcal{NP}$-complete~\cite{wolsey2014integer}. 
Exact techniques to solve this problem usually rely on branch-and-bound (B\&B) algorithms, which explore the discrete variable state space efficiently by pruning nodes associated with infeasible or low-quality solutions~\cite{bertsimasOptimizationIntegers2005a}. 
Various techniques have been developed to improve the performance of B\&B~\cite{morrison2016branch}, including relaxation methods for lower bounds computation~\cite{fisher1981lagrangian} and decomposition into smaller subproblems~\cite{geoffrion1972generalized, barnhart1998branch}. Despite the significant improvements~\cite{achterberg2013mixed}, however, solving MIP problems is still challenging.

The issue is especially problematic in those applications where a MIP problem is solved \emph{online}, i.e., when its solution has to be delivered within a precise time window, as often required to ensure performance, reliability, and safety of a system~\cite{belta2019formal}.
In this context, much work has been devoted to developing schemes to accelerate finding (possibly sub-optimal) solutions, exploiting either some regularity of the problem~\cite{beccuti2004temporal} or domain-specific knowledge~\cite{pippia2019single}. 
Indeed, in many real-time applications, one must repeatedly solve a problem with a fixed structure that depends on external (feedback) information, represented as a set of \emph{parameters} for the problem, as seen in control applications for robotics~\cite{wei2008optimal} or autonomous driving~\cite{quirynen2024real}, and also in state-space models of unknown dynamical systems extracted from data through system identification techniques~\cite{masti2021learning, fabiani2025neural}. In particular, the approach in~\cite{fabiani2025neural} constructs hybrid dynamical models, which are the ones considered in this paper.
This structure can be exploited to partially precompute the solution in an \emph{offline} phase, significantly reducing the operations to be carried out within real-time constraints. 

Within the domain of cyber-physical system control, these approaches are commonly referred to as explicit Model Predictive Control~\cite{alessio2009survey} (MPC). In its original formulation, this technique exploited multi-parametric optimization arguments~\cite{pistikopoulos2020multi} to compute the explicit map between the external input information (the ``parameters'') and the optimal decision associated with them~\cite{alessio2009survey}. 
Explicit MPC, however, suffers from the curse of dimensionality and is therefore impractical for most applications~\cite{bayat2011flexible}, although related works have been devoted to reducing its computational and storage burden by approximating the exact map~\cite{bemporad2006algorithm, malyuta2019approximate}.

A related line of research has focused on constructing a mapping between parameters and \emph{feasible} (though not necessarily optimal) assignments of binary variables~\cite{malyuta2019partition}, allowing numerical solvers to compute the remaining real-valued decision variables during the online phase. These feasible, potentially sub-optimal solutions can be used either to "warm-start" the mixed-integer solvers or fixed directly in the optimization problem. This approach leverages the fact that many optimization problems, such as linear and quadratic programs, can be solved efficiently~\cite{cimini2017exact,nesterov1994interior}.
Numerous studies have explored the use of machine learning to approximate solutions for these types of optimization problems. For example,~\cite{masti2019learning} uses a neural network to learn approximate binary solutions for control tasks, while~\cite{bertsimas2022online} predicts both binary decision variables and the set of active constraints, and~\cite{cauligi2021coco} focuses on predicting the set of relaxed big-M constraints. However, these machine learning approaches often lack a formal soundness guarantee, necessitating fallback strategies to handle cases where the predictor may produce invalid outputs (or \emph{hallucinate}).
In most studies, when provided, fallback strategies rely heavily on domain-specific knowledge. For instance,~\cite{masti2020learning} employs a rule-based approach drawn from established power-engineering principles to handle failures in their ML-based decision-making engine for microgrid applications.

In the context of functional synthesis, we cite Counter-Example Guided Inductive Synthesis (CEGIS), a family of approaches that has achieved significant success due to its ability to certify specification requirements, a critical aspect where machine learning approaches often fall short~\cite{david2017program}.
This iterative, data-driven approach generates functions that meet specified preconditions and postconditions using input-output pairs. CEGIS alternates between synthesis and verification steps: the synthesizer generates a candidate function from a set of template functions using a sample of input-output pairs, while the verifier checks if the function satisfies the required postconditions for all valid inputs. If the function is verified, the process concludes; otherwise, a counterexample that guides the next synthesis step is identified.
CEGIS has been successfully applied across various domains of system verification and control, including the synthesis of controllers~\cite{ravanbakhsh2016robust}, Lyapunov functions~\cite{abate2020formal}, barrier certificates~\cite{peruffo2021automated}, and control Lyapunov functions~\cite{masti2023counter}. 

In this paper, we present a novel application of CEGIS to synthesize an explicit map that links the input parameters of non-convex mixed-integer nonlinear programming problems with polynomial constraints to feasible binary variable assignments. The proposed approach provides a feasibility certificate for the output of the synthesized predictor, addressing the limitations of existing machine learning-based approximations that lack such guarantees.
This map serves as the basis for a suboptimal decision-making engine, significantly reducing the computational load during the online phase.
To our knowledge, this is the first application of CEGIS for generating a feasible binary solution map for MIP problems.

\begin{figure*}[ht]
    \includegraphics[width=0.95\linewidth]{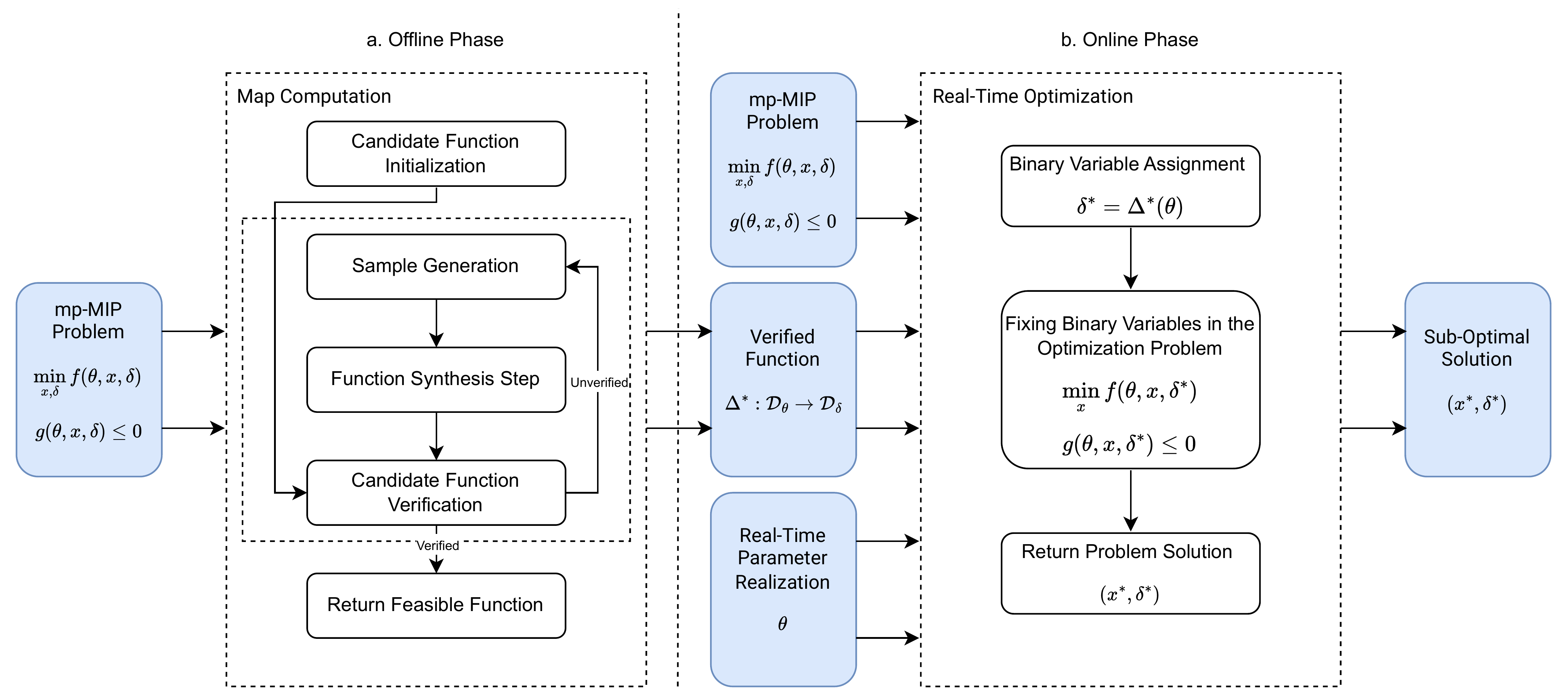}
    \caption{Overview of the proposed two-phase methodology for solving online optimization problems. (a) During the offline phase, a function between real-time parameter realizations and feasible binary variable assignments is synthesized through a CEGIS procedure. This phase is computationally expensive and is executed offline; thus, it does not impact operations during the system execution. (b) In the online phase, the pre-computed map is used to quickly compute feasible values for binary variables, while the remaining continuous decision variables are computed by numerical solvers in real time.}
    \label{fig:overview}
\end{figure*}

The proposed procedure is schematically depicted in \figurename~\ref{fig:overview}. The figure illustrates a two-phase approach for solving mixed-integer programming (MIP) problems, where the computational burden of real-time MIP solving is alleviated by pre-computing part of the problem offline.
In the \textbf{offline phase} (left side of the figure), we initialize the function (or map) $\dfun$ that associates input parameters belonging to a domain to feasible binary variable assignments and we refine it iteratively through cycles of synthesis and verification steps until it is verified.
In the \textbf{online phase} (right side of the figure), the \textit{certified} function $\dfun$ is evaluated for the given parameter values and the binary variables in the MIP problem are fixed to the function output, reducing it to an optimization problem over only continuous variables, for which the numerical solver quickly computes a \textit{sub-optimal solution}.

The remainder of the paper is organized as follows: Section~\ref{sec:preliminaries} introduces the mathematical concepts and notation. In Section~\ref{sec:methodology}, we describe our method. In particular, Section~\ref{sec:offline} describes the offline phase of our procedure, which is meant to construct our binary variables predictor; Section~\ref{sec:online} details how such a function is exploited online. Section~\ref{sec:experiments} evaluates the method using simple numerical examples to assess its practical advantages and limitations compared to techniques for mixed-integer optimization proposed in related works. Section~\ref{sec:limitations} discusses the limitations of the proposed procedure, along with some possible mitigation strategies.

\section{Preliminaries}
\label{sec:preliminaries}

We consider the following multi-parametric mixed-integer nonlinear programming problem (mp-MINLP):
\begin{equation}
    \begin{aligned}
        (x^{\ast}_{\theta}, \delta^{\ast}_{\theta}) := \arg & \min_{x, \delta} && f(\theta, x, \delta) \\
        & \phantom{..} \text{s.t.} && g(\theta, x, \delta) \leq 0, \\
    \end{aligned}
    \label{eq:mpMIP}
\end{equation}
where $\theta \in \Theta \subseteq \mathbb{R}^{N_{\theta}}$, is the parameter vector, $x \in \mathbb{R}^{N_{x}}$ is the continuous decision variable vector, and $\delta \in \{0, 1\}^{N_{\delta}}$ is the binary decision variable vector, with $N_{\theta}, N_{x}, N_{\delta} \in \mathbb{N}$ being the dimensions of corresponding vectors. The function $f$ denotes a nonlinear objective function, while $g$ is a vector-valued polynomial function representing problem constraints, with the inequality evaluated element-wise.
Problem \eqref{eq:mpMIP} is said to be feasible for a given parameter $\theta$ if there exists an assignment of decision variables $(x, \delta)$ that satisfies the constraints $g$. To simplify the exposition, we introduce a predicate, $\feas(\theta)$, representing this condition, defined as follows
\begin{equation*}
    \feas(\theta) \coloneqq \exists x, \delta : g(\theta, x, \delta) \leq 0.
\end{equation*}

Additionally, we introduce the optimization problem in which binary decision variables are fixed to a constant value $\hat{\delta}$. This problem belongs to the class of multi-parametric nonlinear programming (mp-NLP) problems:
\begin{equation}
    \begin{aligned}
        x^{\ast}_{\theta} := \arg & \min_{x} && f(\theta, x, \hat{\delta}) \\
        & \phantom{..} \text{s.t.} && g(\theta, x, \hat{\delta}) \leq 0. \\
    \end{aligned}
    \label{eq:mpMIP_deltafixed}
\end{equation}

Solving this problem using numerical solvers is simpler than solving~\eqref{eq:mpMIP} as it avoids the need to traverse the discrete state space (e.g., via a B\&B strategy~\cite{bertsimasOptimizationIntegers2005a}) to certify feasibility and optimality.

To distinguish between the feasibility conditions of the original mp-MINLP problem~\eqref{eq:mpMIP} and its counterpart~\eqref{eq:mpMIP_deltafixed}, which has fixed binary variable values, we introduce the predicate $\feas_{\hat{\delta}}(\theta)$
\begin{equation*}
    \feas_{\hat{\delta}}(\theta) \coloneqq \exists x : g(\theta, x, \hat{\delta}) \leq 0,    
\end{equation*}
where $\feas_{\hat{\delta}}(\theta)$ evaluates the feasibility of~\eqref{eq:mpMIP_deltafixed} based solely on the existence of at least one assignment of real decision variables $x$ satisfying the constraints $g$.

In what follows, we adopt the assumption~\cite[Assumption~1]{malyuta2019partition}, where the overlap with respect to problem~\eqref{eq:mpMIP} and the set of feasible parameters $\Theta$ is positive, i.e, $\forall \theta \in \Theta$, there exists a neighborhood of the parameter such that its intersection with the feasible region associated with a binary value $\delta$ is non-degenerate.

\subsection{Running example}
\label{ex:example}

The formalism in~\eqref{eq:mpMIP} is commonly used in the context of Model Predictive Control~\cite{garcia1989model} (MPC) to formulate optimal control problems for hybrid dynamical systems. At each time step an optimization problem is solved to compute the optimal decision variables based on the current state of the system and a finite time horizon $N$. The optimal assignment of decision variables for the controls at the first time step is then applied as the action to the system.

To illustrate our approach, we introduce a running example involving a two-dimensional state-space model predictive control problem that tracks a fixed target point. This example features a circular input dead zone and polyhedral obstacle avoidance, both of which are represented using binary variables.
We consider the following two-dimensional piece-wise affine (PWA) controlled dynamical system
\begin{subequations}
    \label{eq:runningExample}
    \begin{equation}
        \begin{bmatrix}
            p_x(k + 1) \\
            p_y(k + 1)
        \end{bmatrix} = 
        \begin{bmatrix}
            p_x(k) \\
            p_y(k)
        \end{bmatrix} 
        + 
        \Delta t \cdot
        \begin{cases}
            \begin{bmatrix}
                0,
                0
            \end{bmatrix}^T & \text{if } \lVert v(k) \rVert_2 < \under{v} 
            \\[1em]
            \begin{bmatrix}
                v_x(k),
                v_y(k)
            \end{bmatrix}^T & \text{otherwise},
        \end{cases}
        \label{eq:dyn_sys_pwa}
    \end{equation}
    where $\Delta t \in \mathbb{R}$ is the discrete time step length, $p(k) = (p_x(k), p_y(k))$ is the system state representing the position in the plane, and $v(k) = (v_x(k), v_y(k))$ is the input velocity, which is known to have a \textit{deadzone}. We impose upper and lower bounds on the position
    \begin{equation}
        \under{p} \leq p(k) \leq \bar{p}~~\forall k,
        \label{eq:state-bounds}
    \end{equation}
    with $\under{p}, \bar{p} \in \mathbb{R}^2$, while we represent the set of possible velocities as 
    \begin{align}
        ||v(k)||_2  \in [\under{v}, \bar{v}] \cup {0}~~\forall k, \label{eq:possible-vel} \
    \end{align}
    where $\under{v}, \bar{v} \in \mathbb{R}$. Additionally, we require this system to remain outside a specified polyhedral \textit{exclusion zone}, thus we impose that 
    \begin{equation}
        \begin{aligned}
            p(k) \notin \mathbf{O}~~\forall k, \\
            \mathbf{O} := \{p \in \mathbb{R}^2 \mid Hp - w \leq 0\},
            \label{eq:poly-exclusion}
        \end{aligned}
    \end{equation}
     where $H \in \mathbb{R}^{m \times 2}$ and $w \in \mathbb{R}^{m}$ represent the $m$ half-spaces that define the polyhedron. 
\end{subequations}

By restricting to a finite temporal horizon (i.e., considering only $k=0,\dots, N$), we transform~\eqref{eq:runningExample} into an equivalent Mixed Logical Dynamical (MLD) representation involving constraints over binary variables which is more suitable for optimization~\cite{bemporad1999control}. This modeling approach introduces auxiliary binary and continuous variables that numerical solvers subsequently optimize. 
We illustrate this modeling procedure by demonstrating how to convert the constraint~\eqref{eq:possible-vel}. Firstly, we impose an upper bound on the velocity norm:
\begin{equation}
    \lVert v(k) \rVert_2 \leq \bar{v}. \label{eq:max-vel}
\end{equation}
We then introduce a binary variable $\delta_u(k)$ at each time step, which indicates whether $\lVert v(k) \rVert_2 < \under{v}$, using the following constraints
\begin{subequations}
    \begin{align}
        M\delta_u(k) &\geq \under{v} - \lVert v(k) \rVert_2 \\
        M(1 - \delta_u(k)) &> \lVert v(k) \rVert_2 - \under{v},
    \end{align}        
    \label{eq:vel-delta-def}
\end{subequations}
where $M \in \mathbb{R}$ is a sufficiently large positive constant used to relax constraints. Notably, when $\delta_u(k)$ evaluates to true, the constraints ensure that $\lVert v(k) \rVert_2 < \under{v}$, while if $\delta_u(k)$ is false, they enforce $\lVert v(k) \rVert_2 \geq \under{v}$. Note that numerical solvers cannot handle strict inequalities, thus they are approximated as non-strict ones in practice.
We can define a real auxiliary variable $z(k) = (z_x(k), z_y(k))$, which is set to $(0, 0)$ when $\delta_u(k)$ is false and to $(v_x(k), v_y(k))$ when $\delta_u(k)$ is true. These conditions are enforced using the following constraints:
\begin{subequations}
    \begin{align}
        -M \delta_u(k) \leq z_x(k) \leq M \delta_u(k) \\
        -M(1 - \delta_u(k)) \leq z_x(k) - v_x(k) \leq M(1 - \delta_u(k)) \\
        -M \delta_u(k) \leq z_y(k) \leq M \delta_u(k) \\
        -M(1 - \delta_u(k)) \leq z_y(k) - v_y(k) \leq M(1 - \delta_u(k)).
    \end{align}
    \label{eq:aux-var-vel-set}
\end{subequations}
We can implement the polyhedral exclusion zone constraints~\eqref{eq:poly-exclusion} by adding a binary variable for each face of each obstacle at each time step. These auxiliary binary variables are evaluated as true if the system state is outside the corresponding plane of the obstacle.
We introduce a binary variable $\delta_o^i$ for each of the $i = 1, \dots, m$ half-spaces defining an obstacle $o$ and enforce their value through the following big-M constraints at each time step $k$:
\begin{subequations}
    \begin{align}
        M \delta_{o}^{i}(k) &\geq H^i p(k) - w^i \\
        M (1 - \delta_{o}^{i}(k)) &> -H^i p(k) + w^i,
        \end{align}        
    \label{eq:obs-delta-def}
\end{subequations}
here $H^i$ represents the row corresponding to the $i$-th half-space of the obstacle polyhedron. When the binary variable $\delta_{o}^{i}(k)$ is true, the condition $H^i p(k) - w^i > 0$ is enforced. Conversely, if it is false, the condition $H^i p(k) - w^i \leq 0$ holds. Additionally, we include a constraint ensuring that at least one of the binary variables associated with the faces of an obstacle must be true
\begin{equation}
    \sum\limits_{i = 1}^{m} \delta_{o}^i(k) \geq 1.
    \label{eq:out-obs}
\end{equation}
Finally, we rewrite the system dynamics in terms of the newly introduced auxiliary variable:
\begin{equation}
    \begin{bmatrix}
        p_x(k + 1) \\
        p_y(k + 1)
    \end{bmatrix} = 
    \begin{bmatrix}
        p_x(k) \\
        p_y(k)
    \end{bmatrix} 
    + 
    \Delta t \cdot
    \left(
    \begin{bmatrix}
        v_x(k) \\
        v_y(k) \\
    \end{bmatrix}
    -
    \begin{bmatrix}
        z_x(k) \\
        z_y(k) \\
    \end{bmatrix}
    \right).
    \label{eq:dyn}
\end{equation}

\paragraph{\textbf{Optimal Control Problem}}

We can now encode our mission as an optimal control problem. We select a quadratic cost function that minimizes both the distance of the state from a target point and the state velocity. In particular, in our implementation, we are tracking the target point $p^{r} = (0, 0)$. The problem objective function is
\begin{equation}
    \sum\limits_{k = 0}^{N - 1} \begin{bmatrix}
        p(k) \\
        v(k)
    \end{bmatrix}^T
    \mathbf{Q}
    \begin{bmatrix}
        p(k) \\
        v(k)
    \end{bmatrix}
    +
    \begin{bmatrix}
        p(N)
    \end{bmatrix}^T
    \mathbf{P}
    \begin{bmatrix}
        p(N)
    \end{bmatrix}.
    \label{eq:obj-func}
\end{equation}
Here, matrices $\mathbf{Q}$ and $\mathbf{P}$ of appropriate dimensions represent weights for the objective function terms. 

The resulting optimization problem belongs to the category of non-convex mixed-integer quadratically constrained problem (MIQCP), which is an instance of the optimization problem~\eqref{eq:mpMIP} and is stated as follows
\begin{subequations}
    \begin{align}
        (x^{\ast}_{\theta}, \delta^{\ast}_{\theta}) := \arg \min_{x, \delta} \quad& ~\eqref{eq:obj-func}
        \\
        \text{s.t.} \quad& 
        ~\eqref{eq:state-bounds} ~\eqref{eq:obs-delta-def} ~\eqref{eq:out-obs} \phantom{.} \forall k = 0..N\label{eq:ex-constr-1}
        \\
        & ~\eqref{eq:max-vel} ~\eqref{eq:vel-delta-def} ~\eqref{eq:aux-var-vel-set} ~\eqref{eq:dyn} \phantom{.} \forall k = 0..(N - 1) \label{eq:ex-constr-2} 
        \\
        & p(0) = \theta,
    \end{align}
    \label{eq:optimization-problem}
\end{subequations}
where $x$ represents the vector of real decision variables, $\delta$ the vector of binary decision variables, and the vector of parameters $\theta$ represents the initial position of the system.

\section{Methodology}
\label{sec:methodology}

In this section, we outline the methodology used in our approach, covering both the offline and online phases. We begin by describing the construction process for the certified predictor, using the optimization problem~\eqref{eq:optimization-problem} as a representative example. Subsequently, we detail the steps of the online phase and discuss the complexity of the transformed optimization problem.

\subsection{Offline Phase}
\label{sec:offline}

\begin{figure*}[ht]
    \centering
    \subfigure[Verification Phase]{
        \includegraphics[scale=0.5]{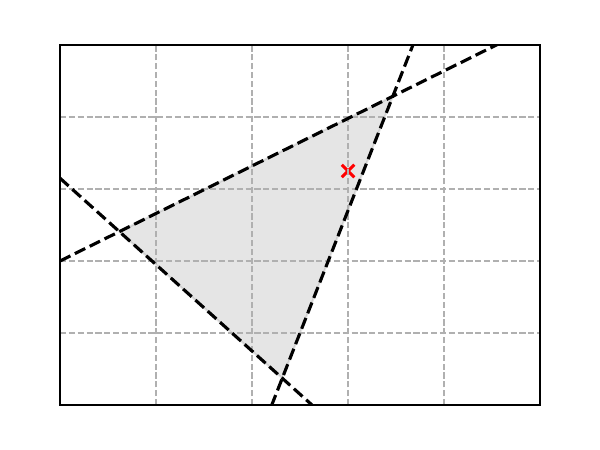}
    }
    \subfigure[Point Sampling]{
        \includegraphics[scale=0.5]{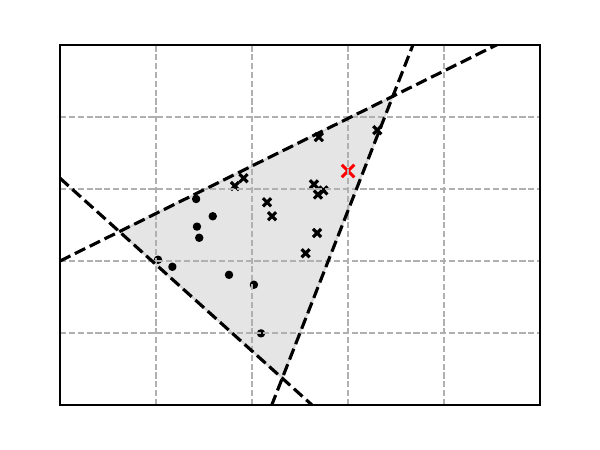}
    }
    \subfigure[Synthesis Step]{
        \includegraphics[scale=0.5]{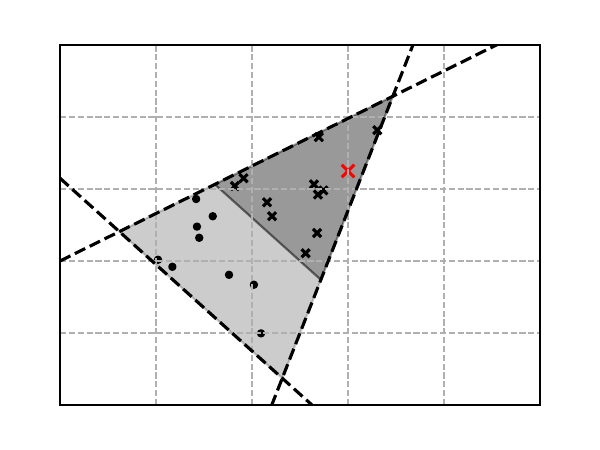}
    }
    \caption{Example iterative verification and refinement process for synthesizing the piecewise function $\dfun$. (a) The verification of the domain $\Omega$ (filled in gray) results in the identification of a counterexample (marked as a big cross). (b) Points from the domain $\Omega$ are randomly sampled for refinement (marked as dots and crosses depending on their optimal binary assignment). (c) The domain $\Omega$ is split into two subdomains $\Omega^{(1)}$ and $\Omega^{(2)}$, forming new branches of the piecewise function.}
    \label{fig:cegis}
\end{figure*}

The offline phase involves constructing the map, as illustrated in \figurename~\ref{fig:overview}(a). 
Our goal is to generate a function $\dfun: \Theta \to \{0, 1\}^{N_{\delta}}$ that accepts parameters within the parameter set $\Theta$ and returns feasible assignments of binary variables for problem \eqref{eq:mpMIP}.
A high-level outline of the algorithm is presented in Algorithm~\ref{al:algorithm-cegis}. The algorithm takes as input a problem instance $P$ having the structure described as in~\eqref{eq:mpMIP}, the parameter set $\Theta$, and an initial candidate function $\dfun$ to be refined. The output is the certified function. In line 3, the $\mathit{refine}$ function uses a learner to synthesize a candidate predictor based on pairs of parameter inputs and their corresponding binary optimizers, obtained by solving the optimization problem~\eqref{eq:mpMIP}. In line 4, the function $\mathit{verify}$ either certifies the correctness of the candidate predictor or identifies a counterexample that violates the desired properties.
\figurename~\ref{fig:cegis} visually illustrates the steps of the iterative process.

\begin{algorithm}[t]
    \SetAlgoNlRelativeSize{0}
    \KwData{$P$, $\Theta$, $\dfun$}
    \KwResult{$\dfun$}
    $\mathit{cex} \gets \mathit{null}$\;
    \Do{$\mathit{cex} \neq \mathit{null}$}{
        $\dfun \gets \mathit{refine}(\dfun, P, cex)$\;
        $\mathit{cex}\gets \mathit{verify}(\dfun, P, \Theta)$\;
    }
    \Return $\dfun$
    \caption{CEGIS procedure}
    \label{al:algorithm-cegis}
\end{algorithm}

\subsubsection{\textbf{Function Structure}}

In this work, we consider functions $\dfun$ that take the form of a piecewise constant map over subsets of the parameter set:
\begin{equation}
    \dfun(\theta) =
    \begin{cases}
        \delta_1 & \theta \in \Omega_1 \\
        & \vdots \\
        \delta_n & \theta \in \Omega_n. \\
    \end{cases}
    \label{eq:pwise-const-func}
\end{equation}
Here, $\Omega_i$, for $i = 1,\ldots,n$, partition the parameter set $\Theta$, that is, $\bigcup_{i = 1}^{n} \Omega_i = \Theta$ and $\Omega_i \cap \Omega_j = \emptyset, \forall i \neq j \in \{1,\dots,n\}$. This template function can be refined indefinitely by further subdividing each subset, allowing for the approximation of functions with arbitrary precision. Various machine learning methods, such as decision trees, k-means clustering algorithms~\cite{macqueen1967some}, or ReLU-activated neural networks~\cite{montufar2014number}, can be used to synthesize an initial candidate predictor in this form.
Moreover, this template function is composed of disjoint regions, allowing the verification and refinement processes to be conducted independently and in parallel at each branch.

Building on the running example introduced earlier (Example~\ref{ex:example}), we illustrate a valid initial candidate function. In particular, a trivial way to initialize the predictor function could involve constructing a piecewise constant function with only one branch, defined over the entire domain $\Theta$, and assigning a random binary variable $\delta$ such that $\dfun(\theta) = \delta, \forall \theta \in \Theta$.

\subsubsection{\textbf{Function Verification}}

The initial candidate function described in the previous section does not provide any formal certification of the feasibility of the returned binary solutions. The critical step in the procedure is the formal verification of the constructed predictor.
We assume that every parameter within the domain $\Theta$ is feasible, i.e, $\feas(\theta)~\forall \theta \in \Theta$. It is important to note that this is not a strict assumption, as we can partition the domain of interest into multiple disjoint subdomains and synthesize a predictor for each.

We frame the requirements of our function in terms of pre-conditions and post-conditions, which are predicates that must be satisfied by the function inputs and outputs, respectively. Specifically, we require that the input parameters are contained within the subset of interest $\theta \in \Theta$. Additionally, we require that each output $\dfun(\theta)$ maintains the feasibility of the optimization problem. This condition can be expressed with the first-order formula $\feas(\theta) \implies \feas_{\dfun(\theta)}(\theta)$. Since we assumed each parameter $\theta$ in the domain is feasible, the implication simplifies to $\feas_{\dfun(\theta)}(\theta)$, which is easier to verify.
Given the desired pre-conditions and post-conditions, we define the first-order formula that must be verified to certify the predictor $\dfun$:
\begin{equation}
    \forall \theta \in \Theta : \feas_{\dfun(\theta)}(\theta).
    \label{eq:validation-formula}
\end{equation}

We are interested in identifying counterexamples where the formula~\eqref{eq:validation-formula} is not satisfied, as these counterexamples enable a more informed refinement process. Therefore, we query the negation of the formula~\eqref{eq:validation-formula}:
\begin{equation}
    \exists \theta \in \Theta : \lnot \feas_{\dfun(\theta)}(\theta).
    \label{eq:neg-validation-formula}
\end{equation}

If the verifier evaluates formula~\eqref{eq:neg-validation-formula} as true, it will return the parameter $\theta^c$  that satisfies the formula, representing a counterexample. Conversely, if the formula is false, we obtain a certificate for the candidate function outputs, allowing us to use our predictor to generate feasible binary variable assignments confidently.

In this work, we focus on certifying piecewise constant functions such as~\eqref{eq:pwise-const-func}, where each region can be verified independently by adjusting the pre-conditions of formula~\eqref{eq:neg-validation-formula} accordingly. Specifically, for a particular branch $i$ of the candidate function, defined in the region $\Omega_i$, we modify the query~\eqref{eq:neg-validation-formula} to search for a counterexample within that specific region:
\begin{equation}
    \exists \theta \in \Omega_i : \lnot \feas_{\delta_i}(\theta).
    \label{eq:neg-validation-formula-branch}
\end{equation}

If the formula~\eqref{eq:neg-validation-formula-branch} is evaluated as unsatisfiable, we can mark that branch as verified. Conversely, if a model satisfying the formula is found, it indicates that the branch domain contains at least one parameter associated with an unfeasible variable assignment and will need further refinement.
Note that formula \eqref{eq:neg-validation-formula-branch} includes polynomial terms involving variables $\theta$ and $x$, and consequently belongs to the theory of \textit{Nonlinear Real Arithmetic} (NRA) with quantifiers. This class of satisfiability problem is decidable using procedures for the quantifier elimination~\cite{tarski1998decision}, ensuring that every query will terminate in a finite number of steps.

Using our running example as a reference, we review the steps involved in verifying a specific branch of the predictor function. Consider the candidate function $\dfun$ introduced earlier, and suppose we want to verify its sole branch. The adapted formula from~\eqref{eq:neg-validation-formula-branch} for verifying the corresponding region is:
\begin{equation*}
    \exists \theta \in \Theta : \forall x : \lnot (\eqref{eq:ex-constr-1} \land \eqref{eq:ex-constr-2}),
\end{equation*}
where the constraints~\eqref{eq:ex-constr-1} and~\eqref{eq:ex-constr-2} are evaluated with fixed $\delta$ values given by $\dfun(\theta)$.

\subsubsection{\textbf{Function Synthesis}}

Suppose that the query of the formula~\eqref{eq:neg-validation-formula-branch} applied to a specific branch $i$ of the piecewise function~\eqref{eq:pwise-const-func} results in a counterexample $\theta^c$. This indicates that the candidate predictor needs to be refined.
We conduct the refinement step using data-driven techniques, starting with the generation of a finite subset of parameters $\hat{\Omega}_i$ from the domain associated with the branch $\Omega_i$. We derive a set of observations $\mathcal{C} = \left\{(\theta, \delta_{\theta}^{\ast}) \mid \theta \in \hat{\Omega}_i\right\}$, where $\delta_{\theta}^{\ast}$ represents the binary solution of the optimization problem~\eqref{eq:mpMIP}.
The synthesizer partitions the unverified region $\Omega_i$ into disjoint regions $\Omega_i^{(1)}, \dots, \Omega_i^{(n)}$. Each newly created region is assigned its corresponding value $\delta_i^{(j)}$ and is added as a new branch in the function~\eqref{eq:pwise-const-func}, replacing the refined branch.

In the context of running example~\ref{ex:example}, assume that during the verification step of the candidate function $\dfun$, a counterexample $\theta^c$ within the region $\Theta$ was found, necessitating the refinement of the function. We start by constructing a set of observations consisting of pairs $(\theta, \delta_{\theta}^{\ast})$ for $\theta \in \Theta$, where $\delta_{\theta}^{\ast}$ denotes the optimal assignment of binary variables obtained by solving problem~\eqref{eq:optimization-problem} for an initial state $\theta$. The synthesizer selects an appropriate hyperplane $(a, b) \in \mathbb{R}^2 \times \mathbb{R}$ to act as a splitter (using, for instance, an iteration of a linear classifier algorithm), thereby dividing the set $\Theta$ into two disjoint sets and refining the candidate function as follows:
\begin{equation*}
    \dfun(\theta) =
    \begin{cases}
    \delta_1 &\text{if } \theta \in \Theta \text{ and } a^{T} \theta - b \leq 0 \\
    \delta_2 &\text{if } \theta \in \Theta \text{ and } a^{T} \theta - b > 0.
    \end{cases}
\end{equation*}

\subsection{Online Phase}
\label{sec:online}

The online phase of the procedure, depicted in \figurename~\ref{fig:overview}(b), is responsible for solving the optimization problem during system execution and closely resembles methods used in related works~\cite{masti2019learning, malyuta2019partition, masti2020learning}. 
Given the parameter values measured at a specific time frame and the map computed during the offline phase, we retrieve a suboptimal feasible binary solution, which is fixed within the original problem~\eqref{eq:optimization-problem}. Numerical solvers are then used to solve the resulting optimization problem for the remaining continuous variables. 

\subsubsection{\textbf{Function Evaluation}}

The first step of the online phase involves evaluating the function $\dfun(\theta)$, which means finding the branch of the piecewise function~\eqref{eq:pwise-const-func} whose domain contains the point $\theta$. 
An efficient representation for encoding the function $\dfun$ is a binary search tree~\cite{tondel2003evaluation}, which allows us to query the function with $\mathcal{O}(h)$ comparisons, where $h$ is the maximum height of the tree. Since each comparison involves evaluating a dot product in $N_{\theta}$ terms, the total complexity of function evaluation is $\mathcal{O}(N_{\theta} h)$.

\subsubsection{\textbf{Online Problem Solution}}

In the subsequent steps of the online phase, the feasible assignment $\delta^{\ast} = \dfun(\theta)$ is used to compute the complete solution to the problem. By fixing the binary variables in problem~\eqref{eq:mpMIP}, we obtain the purely continuous optimization problem~\eqref{eq:mpMIP_deltafixed}, which is guaranteed to be feasible.
Finding the optimal solution for problem~\eqref{eq:mpMIP_deltafixed} in its general form remains an $\mathcal{NP}$-complete problem, as it is a non-convex NLP problem. However, since there is no longer a need to explore the discrete solution space for integer variables, we can focus solely on the continuous solution space using gradient-based methods, which significantly reduce the number of required operations.
In specific cases, such as when problem~\eqref{eq:mpMIP} is a convex Mixed-Integer Quadratic Program (MIQP), the continuous problem~\eqref{eq:mpMIP_deltafixed} simplifies to a convex QP. For this class of problems, efficient solvers suitable for embedded applications exist and can offer worst-case polynomial complexity certification~\cite{cimini2017exact}, effectively reducing the complexity class of the problem during the \textit{online} phase.

\section{Experiments}
\label{sec:experiments}

In this section, we evaluate the effectiveness of the proposed online optimization procedure using numerical examples from hybrid system optimal control as benchmarks.

The primary goal of our proposed procedure is to transform the combinatorial optimization problem~\eqref{eq:mpMIP} into a continuous one~\eqref{eq:mpMIP_deltafixed}. 
To verify the effectiveness of our approach, we measure the reduction in solving time between the original problem and its transformed version.

The synthesized function is designed to approximate the optimal binary solution function, similar to approaches in related works employing machine learning techniques~\cite{masti2019learning, masti2020learning, bertsimas2022online}. However, those methods cannot guarantee the feasibility of the resulting solutions.
To assess this aspect, we first empirically verify that our predictor consistently computes feasible solutions, whereas predictors constructed through machine learning techniques may fail to do so.
Next, we compare the suboptimality of the objective function obtained from solving the continuous problem~\eqref{eq:mpMIP_deltafixed} with integer assignments generated by our predictor against the solution obtained by solving the original problem~\eqref{eq:mpMIP}.

The experiments are conducted on two distinct problem benchmarks. The first benchmark is based on the previously introduced running example, chosen for its inclusion of nonlinear, non-convex constraints that existing methods cannot manage. The second benchmark extends the first, serving as a simple representative example of more commonly encountered problems.

\subsection{Methodology Implementation}
\label{sec:implementation}

We now describe a practical implementation of the algorithm, coded in Python and tested on numerical benchmarks, along with an explanation of the related hyperparameters and their roles. Both the function synthesis procedure and result evaluations were conducted on a machine equipped with an Intel Xeon E7-4830 v4 CPU at 2.0 GHz and 512 GB RAM, running Python 3.12.

The $\dfun$ function is represented as a binary decision tree, initially constructed using an implementation of the \textit{Classification and Regression Tree} (CART) algorithm with the Gini impurity metric, and trained iteratively until an initial specified height $h$ is reached.
To train the initial predictor, we sample $N_i$ points from the parameter space $\Theta$ using the stochastic \textit{Latin Hypercube Sampling} (LHS) technique~\cite{garud2017design} and solve the optimization problem~\eqref{eq:mpMIP} through the Gurobi mixed-integer solver~\cite{gurobi}, treating the resulting binary optimizers as labels for the supervised learning algorithm.
Each leaf node of the constructed decision tree corresponds to a branch of the template function~\eqref{eq:pwise-const-func}, where the predicted value is the most frequent value among the outputs in that node.

To certify the predictor, we use the Z3 SMT solver~\cite{de2008z3} to evaluate whether the formulas~\eqref{eq:neg-validation-formula-branch} have a satisfying model. The verification process takes a leaf of the binary decision tree as input and constructs the corresponding formulas for evaluation. 
We include a small tolerance for formula satisfaction, that is selected to be within the bounds defined by the numerical solvers infeasibility thresholds, ensuring that feasibility guarantees are maintained.
To speed up verification, we run multiple solvers in parallel, with each solver instance operating on the same leaf node. The first decisive result from the portfolio of solvers, whether a satisfying model or a certification of its absence, determines the verification outcome.
The Z3 SMT solver supports quantified formulas in the theory of Nonlinear Real Arithmetic (NRA), enabling it to handle the class of problems discussed in this work.

If a counterexample is identified during the verification step, we solve the problem at that point and generate $N_r$ additional samples: $N_r/2$ points within a small $\varepsilon$-neighborhood around the counterexample, where $\varepsilon$ is adaptively determined as one-sixteenth of the unverified region bounding box length, and $N_r/2$ points from the whole unverified region, using stochastic LHS for both.
We then split the current tree leaf by performing one iteration of the CART algorithm on the unverified branch, continuing until all leaf nodes are verified. Although the process can theoretically refine trees infinitely, for numerical stability, we stop the algorithm if the leaf domain becomes too small or if a configured maximum tree height $h^m$ is reached during synthesis. This does not affect the usability of the technique, as sensors have limited bit resolutions and can only represent a finite set of distinct values.

The algorithm implementation is parallelized to improve efficiency. After the initial synthesis step, each leaf node can be processed independently, as the domains of the leaves are disjoint. We initiate multiple processes that take individual nodes as input, performing verification and synthesis steps as needed. These processes return updated nodes, which are then integrated into the binary tree. The refinement tasks continue until the termination conditions are met.

To further enhance the efficiency of constructing the predictor, as a pre-process step, we divide the parameter domain into multiple subdomains, and we synthesize certified predictors over each one of such domains through the process described above. The predictors are then merged into a single decision tree covering the entire domain of interest. Given the nondeterministic nature of the algorithm, we construct multiple trees for each subdomain and select the best one based on the maximum tree height and the number of nodes.

\subsection{Nonlinear Controller Obstacle Avoidance}
\label{sec:nonlinear-obstacle}

The first benchmark problem involves the optimal control of the dynamical system~\eqref{eq:dyn_sys_pwa} introduced earlier as the running example, with the associated MPC optimization problem~\eqref{eq:optimization-problem}. We consider a time horizon of $N = 3$ time steps, with a time step length of $\Delta t = 1$. The system position is constrained within the region bounded by the lower and upper limits $\under{p} = (0, 0)$ and $\bar{p} = (1, 1)$ at each time step. Additionally, we impose an upper bound on the velocity norm $\bar{v} = 0.08$ and a lower bound $\under{v} = 0.03$, below which the input dead zone is activated. Finally, we include a single rectangular obstacle represented by the lower and upper corner points $o^{(l)} = (0.1, 0.4)$ and $o^{(u)} = (0.7, 0.5)$.

\subsubsection{\textbf{Predictor Construction Details}}

\begin{figure}[t]
    \centering
    \includegraphics[width=.6\linewidth]{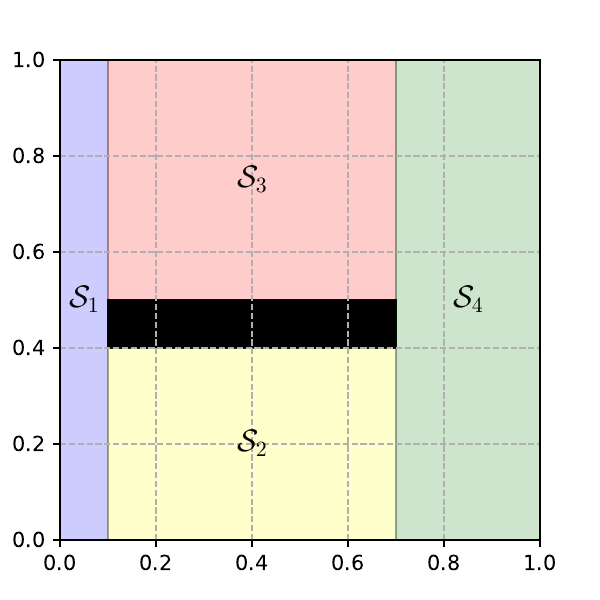}
    \caption{Domain partition for the nonlinear obstacle avoidance problem: the black area represents the infeasible region occupied by the rectangular obstacle, while the remaining states are divided into four feasible regions.}
    \label{fig:obstacle-areas}
\end{figure}

To construct the candidate function, we divide the problem domain into four distinct feasible regions, as shown in \figurename~\ref{fig:obstacle-areas}, and synthesize a predictor for each subdomain. Initially, we sample $N_i = 60$ points to initialize each predictor; during each function refinement phase, we sample an additional $N_r = 20$ points, with half from the neighborhood of the generated counterexample. The decision tree, obtained by merging the best predictors for each subdomain, has a depth of 8 and was constructed using a total of 340 samples.

To establish a baseline for comparison, we compare the synthesized function with traditional machine learning approaches used to approximate the function yielding binary optimizers. Both models were implemented using the scikit-learn library~\cite{pedregosa2011scikit}. Specifically, we trained an uncertified decision tree classifier with a maximum height of 8 (DT) and an uncertified random forest classifier consisting of 10 estimators, each with a maximum height of 8 (RF). Both models were trained using a total of 500 samples. These choices are representatives of predictors lacking formal feasibility guarantees with a comparable online evaluation computational complexity.

\subsubsection{\textbf{Results}}

To evaluate the computational benefits provided by the procedure during online optimization, we sampled 1000 parameters from the feasible domain. We solved both the associated mixed-integer problem~\eqref{eq:mpMIP} (MIP) and the fixed delta problem~\eqref{eq:mpMIP_deltafixed} generated by our certified predictor (CP). Additionally, we considered the continuous problems created by fixing the binary variables to the values predicted by the uncertified decision tree classifier (DT) and the random forest classifier (RF).

\begin{figure}[t]
    \centering
    \includegraphics[width=.7\linewidth,clip,trim={0 0.53cm 0 0 }]{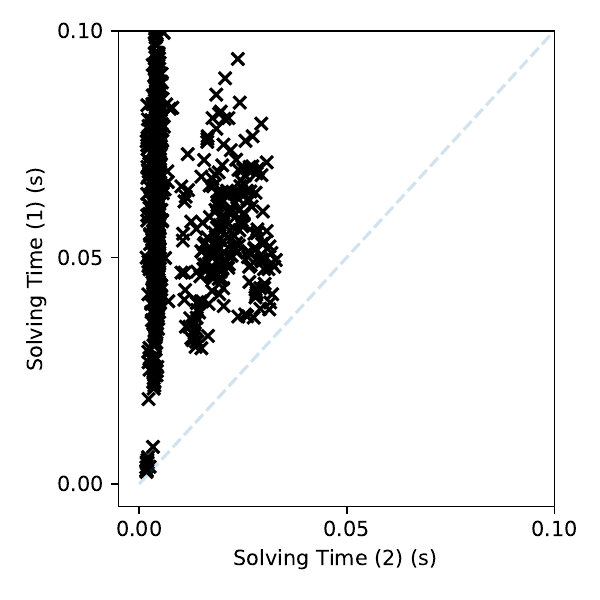}
    \caption{Solving times in the nonlinear system benchmark for optimization problems: \eqref{eq:mpMIP}, on y-axis and \eqref{eq:mpMIP_deltafixed}, with binary variables predicted by $\dfun$ on x-axis.}
    \label{fig:obstacle-times}
\end{figure}

In \figurename~\ref{fig:obstacle-times}, a scatter plot presents the solving times for both the mixed-integer problem and the continuous problem for the same set of parameters. Notably, continuous problems are solved more quickly than mixed-integer ones. This is further confirmed in Table~\ref{tab:obstacle-times-table}, which provides descriptive statistics on the speedup.

\begin{table}[t]
    \centering
    \begin{tabular}{ccc}
        \toprule
        Minimum & Median & Maximum \\
        \midrule
        18.8 \% &  92.4 \% & 97.7 \% \\
        \bottomrule
    \end{tabular}
    \vspace{1ex}
    \caption{Nonlinear system relative decrease in solving times of problem~\eqref{eq:mpMIP_deltafixed} (fixed binary variables) relative to problem~\eqref{eq:mpMIP}.}
    \vspace{-0.5cm}
    \label{tab:obstacle-times-table}
\end{table}

We now empirically verify the soundness of the integer assignments returned by our certified predictor (CP) compared to those returned by models trained using unverified machine learning approaches. Table~\ref{tab:obstacle-feasibility} shows the number of infeasible continuous problems constructed by fixing binary variables to the values predicted by the uncertified classifiers. Notably, the certified predictor (CP) maintains feasibility for all the considered parameters.
Although the number of infeasible continuous problems from traditional machine learning model predictions may seem low, these occurrences can significantly impact online optimization, where the problem is solved multiple times per second.

\begin{table}[t]
    \centering
    \begin{tabular}{ccc}
        \toprule
        \emph{Method} & \emph{\# Feasible} & \emph{\# Unfeasible} \\
        \midrule
        MIP & 1000 & 0 \\
        CP & 1000 & 0 \\
        DT & 978 & 22 \\
        RF & 976 & 24 \\
        \bottomrule
    \end{tabular}
    \vspace{1ex}
    \caption{Nonlinear system feasibility results for continuous problems, comparing the number of feasible instances when fixing binary values predicted by different methods.}
    \label{tab:obstacle-feasibility}
\end{table}

Next, we examine the suboptimality of solutions obtained by solving the continuous problems with fixed binary variables compared to the optimal solutions of the mixed-integer problem~\eqref{eq:mpMIP}. Relative suboptimality is measured using the formula $S^{\theta} = (V^{\ast}_{\text{sub}} - V^{\ast}_{\text{MIP}})/V^{\ast}_{\text{MIP}}$, where $V^{\ast}_{\text{MIP}}$ is the optimal objective value of problem~\eqref{eq:mpMIP}, given that it is considered optimal if the optimality gap with the known lower bound is below $10^{-4}$ and $V^{\ast}_{\text{sub}}$ is the objective value of the corresponding continuous problem~\eqref{eq:mpMIP_deltafixed}. Table~\ref{tab:obstacle-suboptimality} summarizes the results for the various approaches studied, showing the count of suboptimality occurrences within different percentage value ranges for $S^{\theta}$. The first range, labeled "Opt", corresponds to relative suboptimality values less than $10^{-4}$, indicating near-optimal solutions.
Notably, in most cases, the certified predictor $\dfun$ generates assignments of binary variables that result in negligible suboptimality when constructing the continuous optimization problem. Although other tested machine learning approaches slightly outperform CP in terms of solution quality in the considered benchmark, it is important to note that the certified predictor can still undergo further refinement steps while maintaining feasibility.

\begin{table}[t]
    \centering
    \begin{tabular}{ccccc}
        \toprule
        \emph{Method} & \emph{\# Opt} & \emph{\# <1\%} & \emph{\# <5\%} & \emph{\# $\geq$ 5\%} \\
        \midrule
        CP & 820 & 53 & 40 & 87 \\
        DT & 899 & 39 & 28 & 12 \\
        RF & 905 & 33 & 27 & 11 \\
        \bottomrule
    \end{tabular}
    \vspace{1ex}
    \caption{Nonlinear system suboptimality results for continuous problems, comparing the methods by the number of solutions within various suboptimality percentage ranges}
    \label{tab:obstacle-suboptimality}
\end{table}

\subsection{Inertial System with Obstacle Avoidance}
\label{sec:inertial-obstacle}

For the second benchmark, we extend the optimal control problem~\eqref{eq:optimization-problem} previously presented. Specifically, we augment the system state to include the velocity of the object along each axis and control the force applied to the object. To simplify the system, we remove the circular input dead zone while retaining the polyhedral exclusion zone. The system is still tasked with tracking a fixed target point. Consequently, the resulting dynamical system is
\begin{equation}
    \begin{bmatrix}
        p_x(k+1) \\
        p_y(k+1) \\
        v_x(k+1) \\
        v_y(k+1)
    \end{bmatrix} 
    = 
    \begin{bmatrix}
        p_x(k) \\
        p_y(k) \\
        v_x(k) \\
        v_y(k)
    \end{bmatrix} 
    + 
    \Delta t \cdot
    \begin{bmatrix}
        v_x(k) \\
        v_y(k) \\
        a_x(k) \\
        a_y(k)
    \end{bmatrix}
    ,
    \label{eq:inertial-dyn}
\end{equation}
where $\Delta t \in \mathbb{R}$ represents the discrete time step length, $p(k) = (p_x(k), p_y(k))$ denotes the current position in the plane, $v(k) = (v_x(k), v_y(k))$ represents the current velocity, and the newly added term $a(k) = (a_x(k), a_y(k))$ represents the input acceleration. The state position and velocity are bounded:
\begin{equation}
    \under{p} \leq p(k) \leq \bar{p} \quad \under{v} \leq v(k) \leq \bar{v},
    \label{eq:inertial-state-bounds}
\end{equation}
with $\under{p}, \bar{p}, \under{v}, \bar{v} \in \mathbb{R}^2$. Additionally, we impose upper and lower bounds on the acceleration:
\begin{equation}
    \under{a} \leq a(k) \leq \bar{a},
    \label{eq:inertial-input-bounds}
\end{equation}
with $\under{a}, \bar{a} \in \mathbb{R}^2$. The obstacle constraints are encoded using constraints~\eqref{eq:obs-delta-def} and~\eqref{eq:out-obs}, as in the previously introduced non-inertial dynamical system. 

The objective function aims to track the reference point $p^r = (0, 0)$ while minimizing the norm of the velocity and acceleration and is given by
\begin{equation}
    \sum\limits_{k = 0}^{N - 1} \begin{bmatrix}
        p(k) \\
        v(k) \\
        a(k)
    \end{bmatrix}^T
    \mathbf{Q}
    \begin{bmatrix}
        p(k) \\
        v(k) \\
        a(k)
    \end{bmatrix}
    +
    \begin{bmatrix}
        p(N) \\
        v(N)
    \end{bmatrix}^T
    \mathbf{P}
    \begin{bmatrix}
        p(N) \\
        v(N)
    \end{bmatrix}
    .
    \label{eq:inertial-obj-func}
\end{equation}
Here, matrices $\mathbf{Q}$ and $\mathbf{P}$ of appropriate dimensions represent weights for the objective function terms. 

The resulting optimization problem, which falls into the category of Mixed-Integer Quadratic Programming (MIQP), is formulated as follows:
\begin{subequations}
    \begin{align*}
        (x^{\ast}_{\theta}, \delta^{\ast}_{\theta}) := \arg \min_{x, \delta} \quad& ~\eqref{eq:inertial-obj-func} 
        \\
        \text{s.t.} \quad& 
        ~\eqref{eq:inertial-state-bounds} ~\eqref{eq:obs-delta-def} ~\eqref{eq:out-obs} \phantom{.} \forall k = 0..N
        \\
        & ~\eqref{eq:inertial-input-bounds} ~\eqref{eq:inertial-dyn} \phantom{.} \forall k = 0..(N - 1) 
        \\
        & p(0), v(0) = \theta,
        \label{eq:inertial-optimization-problem}
    \end{align*}
    \label{eq:inertial-optimization-problem}
\end{subequations}
where $x$ denotes the vector of real decision variables, $\delta$ represents the vector of binary decision variables, and $\theta$ is the vector of parameters, which corresponds to the initial state of the system, including both the initial position and velocity.

In this model, we consider a time horizon of $N = 5$ time steps, with a time step length $\Delta t = 0.2$. The system position is bounded by lower and upper bounds $\under{p} = (0, 0), \bar{p} = (1, 1)$ at each time step, while velocity within $\under{v} = (-0.05, -0.05), \bar{v} = (0.05, 0.05)$, and the input acceleration within $\under{a} = (-0.05, -0.05), \bar{a} = (0.05, 0.05)$. We consider a single rectangular obstacle, represented by the lower and upper points $o^{(l)} = (0.1, 0.4), o^{(u)} = (0.7, 0.5)$.

\subsubsection{\textbf{Predictor Construction Details}}

We consider a subset of the problem domain $\Theta$ consisting of the rectangle enclosed by the points $\under{\theta} = (0.01, 0.01, -0.01, -0.01)$ and $\bar{\theta} = (0.99, 0.39, 0.01, 0.01)$. We divide the problem domain into the four quadrants of the velocity plane, and synthesize a predictor for each subdomain. Initially, we sample $N_i = 60$ points to initialize each predictor, and during each function refinement phase, we sample an additional $N_r = 20$ points, with half sampled from the neighborhood of the generated counterexample. The decision tree obtained by merging the best predictors for each subdomain has a depth of 10 and was constructed using a total of 900 samples.

For comparison, we also trained an uncertified decision tree classifier (DT) with a maximum height of 10 and an uncertified random forest classifier (RF) composed of 10 estimators, each with a maximum height of 10. Both models were trained using a total of 1000 samples.

\subsubsection{\textbf{Results}}

We conduct these tests similarly to the nonlinear dynamical system case. Specifically, we sample 1000 parameters from the domain used to construct the verified tree and solve both problem~\eqref{eq:mpMIP} (MIP) and the problem~\eqref{eq:mpMIP_deltafixed}. In the latter, integer assignments are generated by our certified predictor (CP), the decision tree classifier (DT), or the random forest classifier (RF). For this test, we use the OSQP solver~\cite{stellato2020osqp} and its extension for solving MIQP problems~\cite{stellato2018embedded}, which are well-suited for embedded settings where optimization problems need to be solved with limited computational resources, as is common in control systems.

\begin{figure}[t]
    \centering
    \includegraphics[width=.7\linewidth,clip,trim={0 0.53cm 0 0 }]{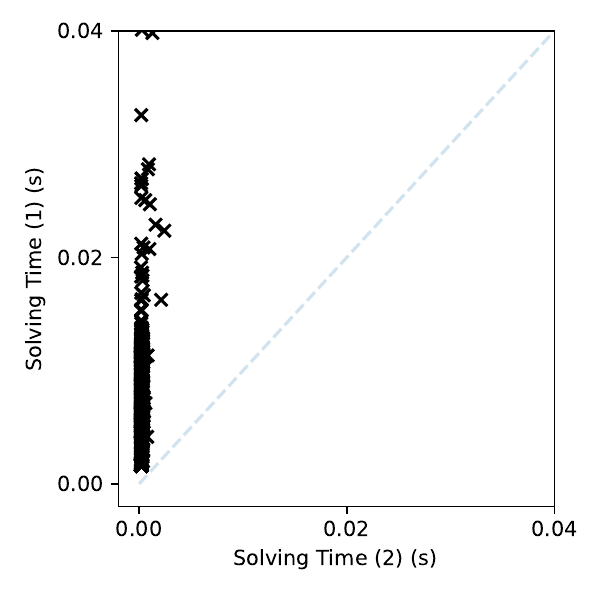}
    \caption{Solving times in the inertial system benchmark for optimization problems: (\eqref{eq:mpMIP}), on y-axis and (\eqref{eq:mpMIP_deltafixed}), with binary variables predicted by $\dfun$ on x-axis.}
    \label{fig:inertial-times}
\end{figure}

In \figurename~\ref{fig:inertial-times}, a scatter plot illustrating the solving times for both the mixed-integer problem~\ref{eq:mpMIP} and the continuous problem~\ref{eq:mpMIP_deltafixed} is presented, while Table~\ref{tab:inertial-times-table} provides the corresponding descriptive statistics. Notably, the reduction in solving times between the two cases is significant, as the exploration of the binary variable state space is eliminated. It is essential to note that, in this case, the problem~\ref{eq:mpMIP_deltafixed} is a convex QP instance, thus finding an optimal solution has polynomial complexity~\cite{cimini2017exact}.

\begin{table}[t]
    \centering
    \begin{tabular}{ccc}
        \toprule
        Minimum & Median & Maximum \\
        \midrule
        81.0 \% & 96.8 \% & 99.7 \% \\
        \bottomrule
    \end{tabular}
    \vspace{1ex}
    \caption{Inertial system relative decrease in solving times of problem~\eqref{eq:mpMIP_deltafixed} (fixed binary variables) relative to problem~\eqref{eq:mpMIP}.}
    \label{tab:inertial-times-table}
\end{table}

We verify the feasibility of integer assignments from our certified predictor (CP) against unverified machine learning models.  Table~\ref{tab:inertial-feasibility} shows the count of infeasible instances for binary variables predicted by the considered classifiers. The certified predictor (CP) maintained feasibility for all parameters, while infeasible predictions occurred with uncertified predictors.

\begin{table}[t]
    \centering
    \begin{tabular}{ccc}
        \toprule
        \emph{Method} & \emph{\# Feasible} & \emph{\# Unfeasible} \\
        \midrule
        MIP & 1000 & 0 \\
        CP & 1000 & 0 \\
        DT & 988 & 12 \\
        RF & 994 & 6 \\
        \bottomrule
    \end{tabular}
    \vspace{1ex}
    \caption{Inertial system feasibility results for continuous problems, comparing the number of feasible instances when fixing binary values predicted by different methods.}
    \label{tab:inertial-feasibility}
\end{table}

We also assess suboptimality by comparing solutions from the continuous problem with fixed binary variables to those from the optimal mixed-integer problem. Table~\ref{tab:inertial-suboptimality} reports suboptimality counts within different percentage value ranges for $S^{\theta}$. Notably, the solution quality from the certified predictor is comparable to that of other uncertified machine learning approaches in this case.

\begin{table}[t]
    \centering
    \begin{tabular}{ccccc}
        \toprule
        \emph{Method} & \emph{\# Opt} & \emph{\# <1\%} & \emph{\# <5\%} & \emph{\# $\geq$5\%} \\
        \midrule
        CP & 979 & 10 & 10 & 1 \\
        DT & 981 & 5 & 2 & 0 \\
        RF & 975 & 9 & 6 & 4 \\
        \bottomrule
    \end{tabular}
    \vspace{1ex}
    \caption{Inertial system suboptimality results for continuous problems, comparing the methods by the number of solutions within various suboptimality percentage ranges}
    \vspace{-0.5cm}
    \label{tab:inertial-suboptimality}
\end{table}

\section{Limitations}
\label{sec:limitations}

During the verification phase, we evaluate the formula~\eqref{eq:neg-validation-formula-branch}, which belongs to the theory of \textit{Nonlinear Real Arithmetic} with quantifiers. This theory is decidable but has doubly exponential time complexity~\cite{collins1976quantifier}. The Z3 SMT solver used in our implementation employs incomplete heuristic decision procedures for the NRA theory with quantifies, which may not always converge. To prevent infinite computations, we impose a timeout on the verification process. In such cases, we mitigate the issue by sampling points from the unverified region and applying the refinement cycle again, although this approach lacks termination guarantees.

Regarding algorithm termination, in~\cite[Definition~4]{malyuta2019partition}, the concept of overlap for an optimization problem~\eqref{eq:mpMIP} and the parameter space $\Theta$ is defined. This condition certifies algorithm divergence if the overlap is zero, indicating that the predictor would require infinite precision. Our procedure cannot synthesize a certified predictor for such problems. One potential mitigation involves detecting these cases and explicitly pre-computing a parameter space partition where each subset has an overlap greater than zero.

Due to the computational complexity of the verification procedure and the curse of dimensionality in constructing the approximator function, our approach is currently suitable only for applications with a limited number of parameters and decision variables. Moreover, at present, constructing a satisfactory predictor requires running the algorithm multiple times and selecting the best one from the results, resulting in a computationally intensive process. Future work could address these limitations by exploring alternative implementations for verifiers and synthesizers, potentially leveraging specific structures inherent in the considered optimization problems. Such improvements would enhance the quality of the predictors, reduce the number of CEGIS iterations required, and ultimately decrease the number of replicas needed to compute a satisfactory predictor.

\section{Conclusion}
\label{sec:conclusion}

We propose an approach based on the CEGIS method to accelerate mixed-integer nonlinear programming problems by computing a map that generates certified feasible integer solutions. The importance of the synthesized function is demonstrated by the ability to construct purely continuous optimization problems, which are easier to solve, albeit possibly leading to suboptimal solutions. We validate the effectiveness of our method through two numerical examples, highlighting the computational time gains in solving optimization problems and empirically verifying both the soundness of the predictor and the solution quality of the returned binary assignments.
Future work will focus on applying post-processing techniques to reduce the number of nodes in the decision tree while maintaining feasibility, thereby improving evaluation times for the function $\dfun$, and refining the predictor to enhance the suboptimality of the returned binary assignments, as well as exploring alternative implementations of verification and synthesis strategies.

\section*{Acknowledgments}

This work has been partially funded by the European Union – NextGenerationEU under the Italian Ministry of University and Research (MUR) through the National Innovation Ecosystem grant ECS00000041 - VITALITY (CUP: D13C21000430001), the SERICS project (SEcurity and RIghts in the CyberSpace, PE0000014, PNRR M4C2 I.1.3), and the 'Resilienza Economica e Digitale' project (CUP: D67G23000060001), awarded as part of the 'Department of Excellence' initiative (Dipartimenti di Eccellenza 2023-2027, Ministerial Decree no. 230/2022). Additional support was provided by project RDS - PTR22-24 P2.1 Cybersecurity. Daniele Masti is also part of INdAM/GNAMPA.

\bibliography{biblio}

\end{document}